\documentclass[10pt,twoside]{article}
\usepackage[a5paper,layoutwidth=148.5mm,layoutheight=7.77817in,
outer=10mm,inner=38pt,tmargin=10mm,bmargin=10mm,footskip=15pt,
]{geometry}

\usepackage[utf8]{inputenc}
\usepackage[lgr,T1]{fontenc}
\usepackage{lmodern}
\usepackage{fix-cm}
\usepackage{textalpha}

\usepackage[comma]{natbib}

\usepackage[french]{babel}
\usepackage[french]{varioref}
\usepackage[autolanguage]{numprint}

\usepackage[hyphens]{url}%
\RequirePackage{etoolbox}
\appto\UrlSpecials{%
  \do\R{\penalty0 \mathchar`\R }%
  \do\U{\penalty0 \mathchar`\U }%
  \do\D{\penalty0 \mathchar`\D }%
  \do\F{\penalty0 \mathchar`\F }%
  \do\A{\penalty0 \mathchar`\A }%
  \do\W{\penalty0 \mathchar`\W }%
  \do\Z{\penalty0 \mathchar`\Z }%
  \do\G{\penalty0 \mathchar`\G }%
  \do\P{\penalty0 \mathchar`\P }%
  \do\N{\penalty0 \mathchar`\N }%
}
\providecommand\url[1]{\texttt{\detokenize{#1}}}

\usepackage{amsmath}

\usepackage{indentfirst}

\ifx\texorpdfstring\undefined\newcommand\texorpdfstring[2]{#1}\fi
\ifx\pdfbookmark[3]{}\fi

\usepackage{booktabs}

\begin{document}
\title{Zénon, les philologues et les scientifiques}
\author{Stefan Neuwirth}
\date{}
\maketitle

Avant de contribuer à l'étude des paradoxes de Zénon, je commence par préciser ma position: celle d'un mathématicien avec une formation philosophique qui enquête sur les circonstances dans lesquelles les mathématiques grecques se sont constituées, et qui découvre un  terrain très controversé dans lequel les sources historiques ne fournissent pas de certitudes.

Cette enquête m'a fait connaitre les quatre fragments que les philologues attribuent à Zénon. 
Trois de ces fragments portent sur la question de la multiplicité: ce qui est, le \og il est\fg, est-il d'un seul tenant, ou peut-on y distinguer une multiplicité de choses qui sont, une multiplicité d'\og étants\fg? Ils ont été moins discutés que les quatre paradoxes du mouvement repris par Aristote dans sa \emph{Physique}: la \emph{dichotomie}, l'\emph{Achille}, la \emph{flèche}, le \emph{stade}.

\section{Introduction. Plaidoyer pour la philologie.}
\label{sec:intr-plaid-pour}

Zénon est un élève de Parménide qui a vécu à Élée sur la mer Tyrrhénienne au 5\ieme\ siècle avant notre ère. Il confronte la pensée du \og il est\fg\ avec les \og opinions des mortels\fg, suscitées par les sensations. Il y a en particulier la sensation d'une multiplicité de choses, qui nourrit l'hypothèse d'une nature discrète de ce qui compose les choses. Zénon mène ses réflexions sous une forme originale qui marque la philosophie au point qu'Aristote le désigne comme le père de la dialectique: il met ces opinions à l'épreuve de la pensée et formule des arguments contradictoires.

Les ouvrages consultés pour cette enquête relèvent davantage de la philologie que de la philosophie: 
j'ai voulu saisir la matière des fragments dans l'épaisseur de leurs mots. Je chercherai à montrer ainsi combien ils sont toujours utiles à la réflexion contemporaine.

Ces fragments reprennent vie dans les analyses des philologues. Leur profession les oblige à rester constamment accrochés à la question du sens et de la place de chaque mot du texte: ils doivent s'assurer de la leçon des manuscrits.

Les philologues sont aussi plus inventifs et aventureux dans leur interprétation des fragments: c'est peut-être parce qu'ils n'ont de compte à rendre qu'au texte lui-même.

Cette présentation a pour ambition une lecture fraiche et directe de textes fondateurs de la philosophie, et de faire revivre les débats enflammés qui ont vu naitre l'axiomatique léguée par les Grecs, voire de comprendre ce qu'il y a d'intemporel dans les questions posées par Zénon.

\section{Comment rendre justice aux paradoxes.}
\label{sec:comm-rendre-compte}

Dans cette section, je présente quelques présupposés épistémologiques qui sont susceptibles de conditionner la lecture des fragments de Zénon proposée dans la section suivante. Ces présupposés se trouvent explicitement dans la \emph{Physique} d'Aristote et concernent la distinction
entre l’acte et la puissance, entre le nombre et la grandeur.

Les quatre fragments originaux de Zénon sont pour nous des moyens d'accéder à l'univers de pensée des éléates. Je propose de les lire au premier degré et de les laisser parler pour eux-mêmes
. Nous en resterons ainsi à la première étape de la stratégie de Zénon: ébranler nos convictions et nous avertir sur le fonctionnement de notre entendement.

Ces fragments s'adressent au fondement d'une théorie de la matière et ne proposent pas une forme pour traiter d'un contenu. En référer à un modèle pour leur répliquer, comme à la théorie des séries infinies, de la mesure ou des nombres réels, c'est se boucher les oreilles et s'empêcher de les écouter.

En fait, les stratégies d'évitement des paradoxes de Zénon consistent à refuser le terrain de l'actuel pour se placer sur le terrain du potentiel, 
c'est-à-dire que le monde n'est pas tel qu'il est, mais qu'il est en devenir. En particulier, la pensée potentielle peut atomiser le monde sans que cette atomisation aboutisse jamais.

Les philosophes et plus encore les mathématiciens ont vite fait d'insérer la pensée de Zénon plus ou moins grossièrement dans un contexte historique qui la désamorce. 
Aristote, quant à lui, résout les paradoxes en soulevant que Zénon confond l'infini du continu --~un segment est infiniment divisible~-- et l'infini des nombres --~il y en a une infinité. Sa définition de l'infini comme un infini en puissance est au fondement de la pensée mathématique de l'infini.

Nous pensons l'infini comme l'abréviation d'un processus inachevé: prenons l'e\-xemple de la \emph{dichotomie}.
\begin{quotation}
  \noindent Le premier argument porte sur l'inexistence du \og se mouvoir\fg, compte tenu du fait que le mobile doit d'abord parvenir à la moitié avant d'atteindre le terme de son trajet\footnote{Aristote, \emph{Physique}, VI~\textsc{ix}~239\textit{b}11-13, voir Jean-Paul Dumont, 1988, page~287.} 
  [\dots].\nocite{dumont88}
\end{quotation}
Nous pouvons penser soustraire d'une chose la moitié de cette chose, puis penser le reste après soustraction de cette moitié, c'est-à-dire l'autre moitié, qui lui est égale. 
Nous pouvons penser la répétition de cette opération. À chaque étape de la 
répétition, nous pouvons penser la chose comme le tout formé de ces moitiés soustraites et du reste après ces soustractions, qui est égal à la dernière moitié soustraite. Mais que pense-t-on quand on prétend penser l'infinité de toutes ces étapes effectuée, c'est-à-dire la chose comme le tout formé de toutes ces moitiés soustraites, et le reste égal à rien?
\begin{itemize}
\item Dit-on seulement ce qu'on a dit juste auparavant, mais en ajoutant la précision que le reste devient arbitrairement petit?
\item Ou dit-on que la chose est vraiment un tout formé d'une infinité de moitiés soustraites?
\end{itemize}
En termes mathématiques, l'égalité suivante, pour tout nombre entier~$n$,
\[1-\frac12-\frac14-\frac18-\dots-\frac1{2^n}=\frac1{2^n}\]
donne lieu à la constatation que pour tout nombre entier~$k$, on a pour tout~$n$ suffisamment grand (il suffit de prendre $n>k$) l'inégalité suivante entre quantités positives: 
\[1-\frac12-\frac14-\frac18-\dots-\frac1{2^n}<\frac1{2^k}\text,\]
ce qu'on abrège et qu'on idéalise depuis les Grecs jusqu'à nos jours par
\[1=\frac12+\frac14+\frac18+\cdots\text.\]

Est-ce que cette égalité a un autre sens que celui d'une abréviation et d'une idéalisation? A-t-on réellement que l'unité est la somme de l'infinité des puissances successives de la moitié? que
\[
\frac1{2^\infty}=0\text{?}
\]
Qu'entend-on par \og réellement\fg?

Voici un autre exemple d'un processus inachevé dont la contrepartie en acte pose question. Nous pouvons penser une chose comme le tout formé de ses $n$ parties $n$ièmes. Mais que pense-t-on quand on prétend que le tout est formé de l'infinité de ses parties infiniment petites?

En termes mathématiques, l'égalité, pour tout~$n$,
\[1=\frac1n+\frac1n+\dots+\frac1n\text{ ($n$ termes)}\]
s'abrège dans le formalisme leibnizien du calcul infinitésimal, en substituant $dx$ à $\frac1n$, par 
\[1=\int_0^1 dx\]
(lire que le tout, l'intervalle [0,1], est l'union de ses portions infinitésimales et que le tout, 1,
est donc la somme des quantités infinitésimales dx lorsque x varie de 0 à 1). Est-ce que cette égalité a un autre sens que celui d'une abréviation? A-t-on réellement que l'unité est la somme de ses quantités infinitésimales?
Notons qu'on ne se permet pas d'abréviation de la forme
\[1=\frac1\infty+\frac1\infty+\cdots\text{ ($\infty$ termes).}\]

La plupart des comptes rendus consistent ainsi en une neutralisation des paradoxes, c'est-à-dire en l'établissement d'un cadre de pensée dans lequel la formulation des paradoxes se dissout en une vacuité conceptuelle.

On peut même voir l'histoire de la réception des paradoxes de Zénon comme l'histoire de propositions de systèmes qui se mettent à l'abri de ces paradoxes. Ainsi Alfred North Whitehead remarque-t-il que \og 
personne n'a jamais effleuré Zénon sans le réfuter, et chaque siècle pense qu'il vaut la peine de le réfuter\footnote{\emph{Symposium in honor of the seventieth birthday of Alfred North Whitehead}, Cambridge, Harvard university press, 1932.\nocite{whitehead32}}\fg{}. Dans ce sens, les paradoxes ont pour objet une limitation de notre entendement à penser les formes de l'espace et du temps.

Du point de vue logique, on aboutit à la conclusion que le système est incompatible avec les paradoxes. Aboutir à cette conclusion est un but légitime pour qui veut se débarrasser des paradoxes pour lui-même ou pour le compte de la théorie qu'il a choisi de défendre.

Mais si on veut rendre justice aux paradoxes, il faut chercher la clé dans le cadre de la pensée de Zénon lui-même. Je propose de lire les paradoxes en amont de toute théorie mathématique et de regarder comment ils conditionnent la mise en place de théories mathématiques. Il faut donc oublier la géométrie euclidienne, suspendre la sommation des séries géométriques, et constater que le texte des paradoxes eux-mêmes ne pose pour seules hypothèses que celles de notre entendement, de notre raison.

Un test qui me parait très utile est de vérifier si nous nous efforçons de penser les paradoxes comme des descriptions de ce qui est et si nous ne les confrontons pas à une image théorique que nous nous sommes faite de ce qui est. Les avertissements suivants par rapport à certaines
théories mathématiques illustrent cet écueil.
\begin{itemize}
\item 
Un intervalle ouvert à gauche~$\mathopen]a,b]$, 
s'il est réellement une chose étante, doit contenir un plus petit élément~$a'$ plus grand que~\(a\) et est donc un intervalle fermé~$[a',b]$. Mais alors comment rendre compte de ce qui sépare $a$~et~$a'$? Voir le début de la section~\ref{sec:les-math-des}.
\item Les nombres réels comme contrepartie mathématique des grandeurs ne sont pas des choses étantes. Un intervalle n'est pas l'ensemble de ses points. Les points géométriques ne sont pas des choses étantes, et même les lignes et les surfaces géométriques n'en sont pas; elles peuvent être pensées comme des limites de choses étantes. Seuls les corps géométriques peuvent prétendre être. Voir la section~\ref{sec:atome}.
\end{itemize}

J'ai donc mis à profit les travaux des philologues pour proposer une traduction plus neutre des fragments qui permet de réfléchir au sens qu'on peut leur donner. J'avertis que le mot grec ἄπειρος, qu'on traduit soit par \og infini\fg{} soit par \og illimité\fg{}, veut littéralement dire \og sans fin\fg{} ou \og sans limite\fg{}; on pourrait le traduire par \og non borné\fg{} pour ne pas conditionner la lecture par une des deux traductions traditionnelles.

Le but de cet article sera atteint s'il permet à ses lecteurs de lire ou de relire les fragments qui suivent d'une manière vivante et ouverte.

\section{Les trois fragments sur la question de la multiplicité.}
\label{sec:le-texte}

Je présente ces fragments dans l'ordre dans lequel Simplicius de Cilicie (490-560) en a rendu compte dans ses \emph{Commentaires sur la Physique d'Aristote}, un millénaire après qu'ils ont été composés par Zénon. On ne sait rien de la transmission de ces fragments jusqu'à Simplicius et ils présentent jusqu'à aujourd'hui de grands problèmes philologiques;
comment comparer le besoin d'authenticité 
entre la fin de l'Antiquité et aujourd'hui?

La traduction que je propose est issue de la comparaison de
traductions existantes \citep{za50,fr55b,un62,vlastos95,ca82,dumont88,simplicius11}.

\begin{quotation}
  \noindent \hfil Fragment~B~2\nopagebreak\smallskip

  Si [un étant qui n'a ni grandeur ni épaisseur ni
  masse]
était ajouté à un autre
  étant, il ne le rendrait en rien plus grand. Car si un étant sans
  grandeur est ajouté à un autre, il n'est pas possible que celui-ci
  gagne en grandeur. Et déjà il s'ensuit que l'ajout égale rien. Puis,
  si la soustraction de quelque chose n'a pas pour effet de rendre
  l'autre chose plus petite, et que son ajout n'a pas pour effet de la
  rendre plus grande, il est clair que l'ajouté comme le retranché
  égalent rien.
\end{quotation}

\begin{quotation}
  \noindent \hfil Fragment~B~3\nopagebreak\smallskip

  Si les étants sont multiples, il est nécessaire que les étants
  soient juste autant qu'ils sont
  , c'est-à-dire ni plus ni
  moins. Mais si les étants sont juste autant qu'ils sont, ils
  devraient être limités.

  Si les étants sont multiples, les étants sont illimités: en effet,
  il y a toujours d'autres étants entre les étants, et à nouveau entre
  ceux-là d'autres encore.  Et ainsi les étants seront illimités.
\end{quotation}

\begin{quotation}
  \noindent \hfil Fragment~B~1\nopagebreak\smallskip

  Si l'étant n'avait pas de grandeur, il ne saurait même pas être.

  Mais s'il est, il est nécessaire que chaque étant ait quelque
  grandeur et quelque épaisseur et que l'un se distingue de
  l'autre. Même raisonnement pour l'étant qui surgit. Celui-ci aura
  aussi une grandeur et une autre partie de lui surgira. Or le dire une fois
  revient à le dire sans cesse. Car aucun de ceux-ci ne sera l'ultime,
  ni l'un ne sera sans relation avec un autre. Donc, si les étants sont
  multiples, il est nécessaire qu'ils soient à la fois petits et
  grands, petits au point de ne pas avoir de grandeur et grands au
  point d'être illimités.
\end{quotation}

\section{Le nombre dans les paradoxes.}
\label{sec:les-math-des}

Zénon nous défie de penser la multiplicité. Qu'est-ce que cela veut dire
d'avoir plusieurs choses devant nous? Comment pouvons-nous arriver au
jugement que les étants sont multiples? Nous devons les imaginer
devant nous comme des lentilles dans une amphore et pouvoir les saisir dans
leur multiplicité en acte, c'est-à-dire d'un coup, dans leur présence. 


Le premier paragraphe du fragment~B~3 nous rappelle que ce qui est, est en acte et non pas en puissance: il nous invite à embrasser les étants dans la pensée, de
sorte qu'il n'en reste pas à ajouter et qu'on n'en ait pas compté de
trop
. Alors les étants sont en nombre fini. Dire qu'il y en a plus qu'il n'y en a n'a pas de sens. Notons que l'hypothèse de la multiplicité n'a pas été utilisée.

En termes mathématiques, Zénon constate que c'est dans la nature du nombre d'être fini. Lorsqu'il exprime que \og si les étants sont multiples, il est nécessaire que les étants soient juste autant qu'ils sont réellement, c'est-à-dire ni plus ni moins\fg, il rappelle dans quelles circonstances un dénombrement fait sens: il faut que les étants à dénombrer soient discernables et donnés avant que l'opération de dénombrement ne commence. Par exemple, une amphore de lentilles contiendra autant de lentilles qu'il y en a, ni plus ni moins.

Le deuxième paragraphe du fragment~B~3 considère l'opération de la pensée qui consiste à distinguer, entre deux étants, un autre étant: \og il y a [\dots] d'autres étants entre les étants\fg. Reportons l'étude de la dimension spatiale de cette opération à la section~\ref{sec:la-topologie-de} et concentrons-nous ici sur le fait qu'une opération a été conçue, ainsi que sa répétition.

Il s'agit ici d'une constatation que je rattache au nombre: chaque nombre a un successeur, et l'opération qui fait passer d'un nombre à son successeur est la version abstraite de toute opération qui peut être répétée. Lorsqu'il mentionne cette opération, Zénon désigne expressément cette répétibilité:
\begin{itemize}
\item \og en effet, il y a toujours d'autres étants entre les étants, et à nouveau entre ceux-là d'autres encore\fg\ (fragment~B~3);
\item \og or le dire une fois revient à le dire sans cesse\fg\ (fragment~B~1).
\end{itemize}
Surtout cette deuxième assertion peut être considérée comme l’essence du principe de récurrence: l'opération qu'on a faite une fois peut être faite sans cesse. Mais notons déjà que cette opération est dite et non faite, c'est-à-dire que c'est une opération que l'esprit conçoit alors que le \og il est\fg\ reste égal à lui-même: il y a une rupture entre la pensée et les étants.

Citons à ce propos un passage célèbre de Poincaré.
\begin{quotation}
  On ne peut donc se soustraire à cette conclusion que la règle du
raisonnement par récurrence est irréductible au principe de
contradiction.

Cette règle ne peut non plus nous venir de l'expérience; ce que
l'expérience pourrait nous apprendre, c'est que la règle est vraie
pour les 10, pour les 100 premiers nombres, par exemple; elle ne
peut atteindre la suite indéfinie des nombres, mais seulement une
portion plus ou moins longue mais toujours limitée de cette suite.

Or, s'il ne s'agissait que de cela, le principe de contradiction
suffirait; il nous permettrait toujours de développer autant de
syllogismes que nous voudrions; c'est seulement quand il s'agit d'en
enfermer une infinité dans une seule formule, c'est seulement devant
l'infini que ce principe échoue; c'est également là que l'expérience
devient impuissante.  Cette règle, inaccessible à la démonstration
analytique et à l'expérience, est le véritable type du jugement syn%
thétique à priori. On ne saurait d'autre part songer à y voir une con%
vention, comme pour quelques-uns des postulats de la géométrie.

Pourquoi donc ce jugement s'impose-t-il à nous avec une irrésis%
tible évidence?  C'est qu'il n'est que l'affirmation de la puissance de
l'esprit qui se sait capable de concevoir la répétition indéfinie d'un
même acte dès que cet acte est une fois possible. L'esprit a de cette
puissance une intuition directe et l'expérience ne peut être pour lui
qu'une occasion de s'en servir et par là d'en prendre conscience\footnote{Henri Poincaré, Sur la nature du raisonnement math{\'e}matique, \emph{Revue de m{\'e}taphysique et de morale}~2, 1894, pages~381-382.}.\nocite{poincare94}
\end{quotation}

La nécessité d'expliciter les articulations d'un raisonnement mathématique comme la récurrence provient de ce que ce dont on traite n'est pas accessible à la sensibilité: il s'agit d'une réflexion suggérée par l'objet de l'étude, les choses étantes.

Je constate à ce niveau le contenu mathématique suivant des paradoxes de Zénon sur la multiplicité.
\begin{itemize}
\item Si un étant ou une multiplicité d'étants admet une opération qui fait augmenter le nombre d'étants, on est en contradiction avec l'assertion que les étants sont juste autant qu'ils sont. En d'autres mots, on ne peut pas penser ce qui est donné comme consubstantiel à une opération qui l'augmente, qui l'altère.
\item La pensée d'une opération qui peut être répétée n'a pas de borne et cette répétition s'étend aussitôt en une succession illimitée de ces opérations. Cela fait apparaitre au grand jour la rupture entre la pensée en puissance et les étants en acte.
\end{itemize}

\section{La grandeur dans les paradoxes.}
\label{sec:atome}

Le premier paragraphe du fragment~B~1 rappelle qu'un étant participe à l'être seulement dans la mesure où il participe à la grandeur. Il s'agit d'expliquer l'étendue, l'extension, la grandeur des choses étantes.

Le fragment~B~2 cherche à penser un étant sans grandeur ni épaisseur ni
masse. C'est un étant hypothétique qui aurait les propriétés du point
géométrique. Or la grandeur, l'épaisseur, la masse sont des propriétés
communes aux objets que nous percevons. Zénon demande quelle serait la contribution d'un étant sans grandeur à un étant qui a quelque grandeur: 
ajouté ou enlevé, un étant sans grandeur ne pourrait altérer
l'autre parce qu'il ne peut pas contribuer à sa grandeur; il est donc littéralement rien.

La conclusion de ce raisonnement pose de grandes difficultés. Elles sont très bien décrites par Aristote.
\begin{quotation}
  Admettons que [ce corps] soit divisé. Que sera le reste? une grandeur? Impossible, car il y aurait quelque chose de non divisé alors que le corps était totalement divisible. Or, si le reste ne consiste ni en un corps ni en une grandeur, et que pour autant la division existe, ou bien le corps sera formé de points, et ses constituants seront des éléments dépourvus de grandeur, ou bien ce ne sera rien du tout; ce qui reviendrait à dire qu'il aurait été engendré et se trouverait composé à partir du néant, et que le tout ne serait rien d'autre qu'un fantôme phénoménal\footnote{Aristote, \emph{De la génération et de la corruption}, I, \textsc{ii}, 316~\emph{a}~13, voir Jean-Paul Dumont, 1988, page~773.}.
\end{quotation}
On peut même se demander si Zénon ne se fait pas l'avocat d'une analyse logique pure qui mène à l'impossibilité de toute image qu'on serait tenté de se faire: il propose un raisonnement limpide qui s'applique à toute image potentielle d'un corps divisible et l'abolit. Je vais néanmoins essayer de décrire de telles images spatiales sur la base du texte de Zénon.

La première image qu'on voudrait pouvoir se faire est celle d'un étant sans grandeur. Lorsqu'on casse une mine de tungstène, la pointe en est un unique atome de tungstène: c'est un étant très petit, et pourtant sa taille limite la résolution des microscopes à effet tunnel: il a une certaine grandeur, de l'ordre de l'ångström ($10^{-10}$~m). Une autre image pourrait être celle de la particule élémentaire, dont je peux seulement constater que la question de son extension spatiale me laisse perplexe.

Les mathématiciens proposent des objets mathématiques qui n'ont \og ni grandeur ni épaisseur ni masse\fg: les surfaces, les lignes, les points. Il s'agit de concepts qui sont l'idéalisation de limites de solides, de limites de surfaces, de limites de lignes. Le fragment~B~2 montre que même si notre entendement use de ces idéalisations avec une grande facilité, elles ne participent certainement pas aux étants: elles n'ont pas de grandeur et si on les retranche en pensée à un étant, on retranche littéralement rien.

La théorie mathématique actuelle qui cherche à répondre à ces questions est la théorie de la mesure. Pour elle aussi, un ensemble a une grandeur seulement si ses constituants en ont une. On peut chercher à constituer un objet à partir de très petits objets, comme des tranches, mais penser ces tranches comme des plans comme le fait la méthode des indivisibles de Cavalieri est une idéalisation qui n'a pas lieu en acte.

\section{Vérité et opinion chez Parménide.}
\label{sec:la-philosophie-des}

Le témoignage d'Aristote dans la \emph{Métaphysique} que \og c'est donc à bon droit que celui qui, le premier, trouva un art quelconque, dégagé des sensations communes, excita l'admiration des hommes\footnote{Aristote, \emph{Métaphysique}, 981\emph{b}13-15, voir Jules Tricot, 1953, page~8.}\fg\nocite{ar86} nous pousse à rechercher l'origine des mathématiques grecques dans des préoccupations métaphysiques et en particulier à voir un lien entre la méfiance vis-à-vis des sensations et le projet de théories déductives.

Les historiens de la philosophie font remonter ces préoccupations au
poème de Parménide, né à Élée à la fin du 6\ieme\ siècle avant notre
ère. Ce poème est considéré comme l'œuvre inaugurale de la métaphysique. Il
raconte une initiation. 


Après un long voyage qu'il décrit minutieusement, Parménide rencontre
la déesse.
\begin{quotation}
\og Jeune homme, toi qui viens ici, accompagné

De cochers immortels, portés par des cavales,

Salut! [\dots] Apprends donc toutes choses,

Et aussi bien le cœur exempt de tremblement

Propre à la vérité bellement circulaire,                        

Que les opinions des mortels, dans lesquelles

Il n'est rien qui soit vrai ni digne de crédit;

Mais cependant aussi j'aurai soin de t'apprendre

Comment il conviendrait que soient, quant à leur être,

En toute vraisemblance, lesdites opinions, 

Qui toutes vont passant toujours\footnote{Parménide, fragment~B~1~28-32, voir Jean-Paul Dumont, 1988, page~256.}.
\end{quotation}

La vérité se dit ἀλήθεια 
en grec et signifie
\og dévoilement\fg\ au sens propre; l'opinion se dit δόξα
, mot qui a donné le suffixe \og -doxe\fg\ en français. La déesse oppose donc
l'une à l'autre: l'une est \og exempte de tremblement\fg, les autres \og vont
passant toujours\fg.  L'opinion se base sur les sensations et
l'apparence; pour accéder à la vérité, à ce qui est (dans un sens
absolu), il faut chercher à le penser. La déesse propose de les enseigner
toutes les deux, mais les quatre derniers vers sont surprenants: les
opinions peuvent être vraisemblables et l'homme peut les organiser de
manière à les rendre probables; la déesse lui enseignera aussi la
manière dont les opinions peuvent prendre cette forme convenable. Cela
reste une nécessité pour l'homme qui, même si on lui a dévoilé la
vérité, continuera à user de ses cinq sens et voudra donner de la
vraisemblance à ses sensations.

D'un autre côté, les opinions peuvent seulement prétendre à la vraisemblance et elles sont changeantes.


La déesse présente ainsi la
vérité dévoilée.
\begin{quotation}
Mais il ne reste plus à présent qu'une voie

Dont on puisse parler: c'est celle du \og il est\footnote{Parménide, fragment~B~8~1-2, voir Jean-Paul Dumont, 1988, page~261.}\fg.
\end{quotation}

Guido Calogero propose d'accéder à la leçon de Parménide par une
analyse logique des assertions: \og [la pensée] ne peut s'effectuer sans
prendre corps dans une affirmation, dont le verbe être est le
constituant essentiel, parce qu'il est la forme commune et universelle
de toute qualification ou prédicat possible, qui ne peut être pensé
sans celle-ci, alors que celle-ci peut réciproquement être pensée pour
elle-même et sans aucun ajout\footnote{Guido Calogero, \emph{Studi sull'eleatismo}, Florence, La nuova Italia, 1977, pages~8-9, ma traduction.}\nocite{ca77}\fg.  Je suis pour
ma part la voie de Parménide en concentrant ma pensée sur le fait même
d'être en faisant abstraction des étants. 

Reportons l'étude de la voie du \og il est\fg\ à la section suivante et lisons comment la déesse introduit aux opinions des mortels.

\begin{quotation}
\indent Mais ici je mets fin au discours assuré

Ainsi qu'à la pensée visant la vérité.

Désormais apprends donc l'opinion des mortels

En ouvrant ton oreille à l'ordre harmonieux

Du discours composé pour ton enchantement\footnote{Parménide, fragment~B~8~50-52, voir Jean-Paul Dumont, 1988, page~263.}.
\end{quotation}
Comme promis, la déesse enseigne à présent l'opinion, le monde des
sensations. L'opposition par Parménide de la vérité aux opinions
continue à enflammer les philologues et les philosophes. 
Ce n'est pas un discours de vérité,
mais ce discours est nécessaire pour rendre compte de la nature, des
sensations, du temps. 

Mario Untersteiner 
soutient
que le discours de l'opinion qui suit le fragment~B~8 n'est
pas l'objet de la condamnation de Parménide: être et opinion sont sur
une même voie qui mène l'être à la vérité et l'opinion à la
vraisemblance. La vérité est définitive; la vraisemblance par contre
est un chemin pour toujours inachevé.

Il replace les fragments de Zénon dans le contexte de
la pensée de Parménide. Selon lui, les paradoxes s'appliquent à la
δόξα (l'opinion) et non à l'ἀλήθεια (la vérité dévoilée): les opinions ne se décrivent pas rationnellement à l'aide de la multiplicité\footnote{Mario Untersteiner, \emph{Zenone. Testimonianze e frammenti}, Florence, La nuova Italia, 1962, pages~186-187.}\nocite{un62}.

\section{La description du \og il est\fg\ chez Parménide.}
\label{sec:la-physique-des}



Pour parler de la voie du \og il est\fg, la déesse en fournit une description à laquelle nous accédons par notre intuition de l'espace.
\begin{quotation}
Sur cette voie il est de fort nombreux repères,

Indiquant qu'échappant à la génération,

Il est en même temps exempt de destruction:

Car il est justement formé tout d'une pièce,

Exempt de tremblement et dépourvu de fin.

Et jamais il ne fut, et jamais ne sera,

Puisque au présent il est, tout entier à la fois,

Un et un continu. Car comment pourrait-on

Origine quelconque assigner au \og il est\fg?

Comment s'accroitrait-il et d'où s'accroitrait-il\footnote{Parménide, fragment~B~8~2-7, voir Jean-Paul Dumont, 1988, page~261.}?
\end{quotation}
On accède au \og il est\fg\ par une méditation et celle-ci a lieu dans une étendue qui permet au \og il est\fg\ d'être pensé comme \og un et un continu\fg. 
Nous nous concentrerons sur le passage suivant qui étudie 
sa topologie.

\begin{quotation}
\leavevmode\hphantom{Apparait impossible. }Et il n'est pas non plus

Divisible en effet, puisqu'il est en entier,

Sans avoir çà ou là quelconque chose en plus

Qui pourrait s'opposer à sa cohésion,

Ou quelque chose en moins. Il est tout rempli d'être.

Aussi est-il tout continu. En effet, l'être

Embrasse au plus près l'être\footnote{Parménide, fragment~B~8~22-25, voir Jean-Paul Dumont, 1988, page~262.}. 
\end{quotation}
Reprenons les termes qui décrivent le \og il est\fg\ dans son extension. Il
est \og tout rempli d'être\fg, \og tout entier à la fois\fg, \og tout continu\fg,
\og formé tout d'une pièce\fg, rien ne \og pourrait s'opposer à sa
cohésion\fg. Il n'est \og pas non plus divisible\fg. Ernst Cassirer constate qu'\og en vertu de cette conclusion, le concept mathématique de l'espace
``vide'', tel qu'il a été pensé par les pythagoriciens comme
``écartement'' (διάστημα) entre les corps, est dévalué ontologiquement et aboli\footnote{Ernst Cassirer, \emph{Die Philosophie der Griechen von den Anf\"angen bis Platon}, Berlin, Ullstein, 1925, page~42, ma traduction.}\fg\nocite{ca25}.
\begin{quotation}
Car la Nécessité puissante le retient

Dans les liens l'enchainant à sa propre limite;

C'est pourquoi la sentence a été décrétée

Que l'être ne saurait se dispenser de fin. 

Il ne manque de rien, en effet
\footnote{Parménide, fragment~B~8~30-33, voir Jean-Paul Dumont, 1988, page~262.}.
\end{quotation}

\begin{quotation}
\indent Mais puisque existe aussi une limite extrême,

Il est de toutes parts borné et achevé\footnote{Parménide, fragment~B~8~42-43, voir Jean-Paul Dumont, 1988, page~262.} [\dots]


\end{quotation}




La première leçon de cette expérience est que le néant, le non-être ne peut pas être pensé: il n'est donc pas. Donc l'être est formé d'un seul tenant. 
%
De même,
\begin{itemize}
\item il échappe à la génération et à la destruction; il est immobile;
\item il est exempt de tremblement, un continu indivisible, gonflé à l'instar d'une balle bien ronde;
\item il est sans commencement et sans fin;
\item il est de toutes parts borné et achevé avec une limite extrême.
\end{itemize}

Le \og il est\fg\ fait ainsi émerger 
des concepts topologiques: 
les oppositions continu\slash discret, à
bord\slash sans bord, limité\slash illimité, connexe\slash formé de plusieurs parties. 

\section{La description de la multiplicité chez Zénon.}
\label{sec:la-topologie-de}

Dans le deuxième paragraphe du fragment~B~3, on cherche à concevoir qu'il y a au moins deux étants. Comment va-t-on les concevoir? Nous dirions aujourd'hui qu'on les conçoit spatialement, selon leur extension. Mais alors cela veut dire que nous les concevons séparés, que l'un est ici, et l'autre là. 
Cet écartement correspond à un troisième étant et nous sommes ainsi obligés de penser un étant supplémentaire.

Voici comment Hermann \citet{fr55b} analyse ce fragment.
\begin{quotation}
  Selon Zeller [\emph{La philosophie des Grecs considérée dans son développement historique}, volume~2, page~75] et d'autres, Zénon veut dire: deux choses ne pourraient être deux que si l'une est séparée de l'autre, et s'il faut qu'elles soient séparées, quelque troisième chose doit se trouver entre elles. Autant que je voie, il n'y a pas de support pour une telle explication et elle ne parait pas probable non plus. Nulle part les éléates nient que deux choses, s'il y avait multiplicité, ne pourraient être en contact. Et si Zénon parle de multiplicité, il inclut toujours la divisibilité; mais par divisibilité il n'entend pas seulement la séparabilité physique, mais aussi la possibilité d'une distinction mentale entre parties ou domaines à l'intérieur de l'objet d'un seul tenant. Il n'y a donc littéralement pas d'espace pour l'insertion d'un diaphragme physique qui tiendrait séparés deux domaines limitrophes; et il n'existe à fortiori pas de raison logique de postuler qu'une troisième chose s'interpose quand il faut en différencier deux. \emph{Res rem finire videtur [l'objet limite l'objet
    ]} (Lucrèce [De la nature] I~998).

  Pour expliquer la deuxième partie du fragment, il ne reste donc que l'autre alternative, c'est-à-dire d'opérer avec des choses (ou des parties) de grandeur indéfinie. Si nous attribuons à quelque unité multiplicité, c.-à-d.\ divisibilité, une certaine partie de celle-ci est ici et une autre là. Mais même si l'ici et le là sont très proches l'un de l'autre, nous sommes libres de poser le quelque chose ici et le quelque chose là si petits qu'on peut encore coincer un troisième entre les deux. Cette opération peut être répétée un nombre quelconque de fois sans qu'une limite soit jamais atteinte. Car la prémisse de la multiplicité et de la divisibilité (mentale) ne permet pas l'hypothèse d'une ultime unité indivisible, et une logique stricte n'autorise pas une transition progressive du très petit vers le sans extension. Car c'est le point essentiel de beaucoup d'arguments de Zénon, que l'absence de grandeur doit être strictement différenciée de toute grandeur, aussi petite soit-elle\footnote{Hermann Fränkel, \emph{Zenon von Elea im Kampf gegen die Idee der Vielheit}, dans \emph{Wege und Formen fr\"uhgriechischen Denkens}, Munich, C. H. Beck, 1955, pages~201-202, ma traduction.}.\nocite{fr55b}
\end{quotation}

Le deuxième paragraphe du fragment~B~3 approfondit la pensée de la multiplicité. On peut aussi l'interpréter ainsi: à partir du moment où on pense que l'être est divisé en étants, cette divisibilité de l'être s'appliquera de manière homogène à chaque étant. Donc la pensée de la divisibilité désintègre immédiatement l'être en une infinité d'étants!

Les fragments décrivent donc deux mouvements d'analyse de l'étant: des parties données vers le tout et du tout vers ses parties. Ils démontrent que chacun de ces mouvements d'analyse est infini.


Le deuxième paragraphe du fragment~B~1 est lui aussi sujet à deux interprétations contradictoires qui mènent à deux traductions différentes.
\begin{itemize}
\item Les étants se distinguent et on peut penser ce qui est
  à côté d'un étant. Mais alors les étants s'étendent à l'infini et la
  multiplicité d'étants est illimitée.
\item On peut y voir une description matérielle de la deuxième interprétation du fragment~B~3. Dans ce cas, chaque étant est pensé comme composé de deux parties distinctes qui se côtoient et le \og même raisonnement\fg\ est appliqué à une des deux parties. On aboutit à un nombre illimité de parties. On a à nouveau une descente infinie; dans cette interprétation,  les deuxièmes paragraphes des fragments~B~1 et~B~3 se rejoignent.
\end{itemize}

Hermann \citet[pages~223-225]{fr55b} propose une interprétation étonnante du fragment~B~1: selon lui,
Zénon considère successivement la surface de l'étant, la surface de sa
surface, etc.\ comme des pellicules enchâssées et montre ainsi que l'étant n'a pas de limite ultime.
\begin{quotation}
  Comme l'épaisseur est nommée, l'écartement réfère ici à l'extension en profondeur de l'objet, et les parties en contraste seront les faces opposées ou surfaces du corps dans la direction du regard: devant et derrière. L'épaisseur suppose, dit Zénon, un certain écartement entre l'avant et l'arrière. Cette explication est confirmée par la parallèle dans la deuxième phrase. Zénon ne parle certes ici que de μέγεθος [grandeur] seule, sans ajouter καὶ πάχος [et profondeur] pour spécifier la dimension; mais à la place du ἀπέχει \og est à l'écart de\fg{} neutre, il utilise cette fois-ci προέχει \og s'élance, est situé plus à l'avant\fg{} -- une expression qui réfère spécialement à l'extension et à l'écartement en profondeur.

  De ce parallélisme résulte de plus que αὐτοῦ [même] a la même fonction dans la deuxième version que dans la première et la troisième. Il ne se rattache pas à προέξει (\og quelque chose sera plus à l'avant que ceci\fg{}), mais à τι et est partitif-possessif: \og quelque chose de lui (c.-à-d.\ du προέχον [se situant plus à l'avant]) se situera plus à l'avant (qu'une autre partie du προέχον)\fg{}. La face avant de l'objet a donc elle-même de nouveau une face avant ou surface, qui se situe plus à l'avant que la face arrière de la face avant de l'objet. Selon la prémisse que tout ce qui existe a trois dimensions, la surface d'un corps, si elle existe, doit être comprise comme une couche ou peau, et a elle-même à nouveau une surface ou peau. Bien sûr que c'est ridicule; mais nous ne devons pas oublier que c'était justement l'intention de Zénon de ridiculiser la théorie de la multiplicité\footnote{\emph{Ibid.}, pages~223-225, ma traduction.}.
\end{quotation}

Ces fragments posent davantage de questions qu'ils ne suggèrent de réponses! Le lieu est un étant comme un autre, contrairement à notre conception de l'espace comme un réceptacle neutre.

Voici un tableau de traductions de trois mots clés des fragments de Zénon.\bigskip

\centerline{\footnotesize\setlength\tabcolsep{3pt}%
  \begin{tabular}{@{}rlll@{}}
    \toprule
    &ὄντα&ἀπέχειν&προέξει\\
    &étants&se distinguer&surgir\\\midrule
    \citet{fairbanks98}&things&to be separated&to precede\\
    \citet{Lee36}&things&to be at a certain distance&to be beyond\\
    \citet{za50}&choses&se trouver à une certaine distance&être situé en avant\\
    \citet{dielskranz51}&seiende Dinge&Abstand haben&davor liegen\\
    \citet{fr55b}&seiende Dinge&Abstand haben&nach vorne liegen\\
    \citet{un62}&realtà&distinguersi&estendersi accanto\\
    Vlastos (1971)&existents&to extend beyond&to project\\
    \citet{ca82}&étants&être en-dehors de&précéder\\
    Dumont (1988)&existants&avoir une certaine distance&être devant lui\\
    \citet{verelst06}&beings&to keep away&to jut out\\
    \citet{koehler14}&Dinge&an anderer Stelle liegen&weiter vorn liegen\\
    \citet{laksmost16}&things&to be distinct&to be more\\
    \citet{laksmost16}&choses&être distinct de&être en surplus\\\bottomrule
  \end{tabular}}\bigskip

Et voici les tentatives de traduction de la préposition πρὸς dans un syntagme dont Fränkel affirme qu'il est clair, mais qu'il ne se laisse pas traduire de manière satisfaisante\footnote{\emph{Ibid.}, page~225, note~4.}. Aucun traducteur ne tente d'exploiter un autre sens de cette préposition, le sens spatial de \og vers\fg.\bigskip

\centerline{\footnotesize\setlength\tabcolsep{3pt}%
  \begin{tabular}{@{}rl@{}}
    \toprule
    &ἕτερον πρὸς ἕτερον οὐκ ἔσται\\
    &ni l'un ne sera sans relation avec un autre\\\midrule
    \citet{fairbanks98}&nor will one thing differ from another\\
    \citet{Lee36}&nor out of relation to another\\
    \citet{za50}&et il n'arrivera jamais qu'on ne puisse opposer une de ces parties à\\
    &une autre\\
    \citet{dielskranz51}&nie wird der eine ohne Beziehung zum anderen sein\\
    \citet{fr55b}&noch wird es dahin kommen, daß eine Gegenübersetzung zweier Dinge\\
    &unmöglich ist\\
    \citet{un62}&né una di esse non avrà relazione con l'altra\\
    Vlastos (1971)&or without one [part] similarly related to the other [part]\\
    \citet{ca82}&ni telle qu'il n'y ait pas de relation d'une de ses parties à l'autre\\
    Dumont (1988)&et il n'est aucun existant qui n'existe pas en relation avec un autre\\
    \citet{verelst06}&nor the one will be unrelated to another one\\
    \citet{koehler14}&noch wird eines nicht in Beziehung zu einem anderen stehen\\
    \citet{laksmost16}&nor will there be any part of it that will not be in relation with another\\
    \citet{laksmost16}&et il n’y aura rien qui ne soit en relation avec un autre\\\bottomrule
  \end{tabular}
}

\section{Les constats incompatibles.}
\label{sec:les-const-incomp}

Notre lecture de Zénon donne que penser la chose étante, en acte, aboutit à ce qu'elle
\begin{itemize}
\item n'est pas discrète parce que sinon la chose étante est atomisée en indivisibles sans grandeur,
\item et n'est même pas formée de plusieurs étants. Elle est \emph{une} et \emph{d'un seul tenant}.
\end{itemize}

Les théories physiques et mathématiques qui passent outre deviennent donc irréalistes sur un de ces trois points au moins:
\begin{itemize}
\item les étants sont infinis, en acte;
\item des individus sans grandeur constituent de la grandeur;
\item elles décrivent le réel comme un devenir, en puissance et non en acte.
\end{itemize}

\section{Conclusion. La place des paradoxes dans l'histoire des mathématiques.}
\label{sec:la-place-des}

La rareté et la forme sèche des fragments présocratiques discutés ici n'empêchent pas de saisir le continu et l'illimité comme défis majeurs pour les éléates et les pythagoriciens. Les paradoxes de Zénon ont décrit très tôt les obstacles d'une description vraisemblable du rapport entre le continu et ses constituants. 

Proposer un exposé sur les paradoxes de Zénon à une conférence d'histoire des mathématiques est cependant un exercice risqué. Il doit naviguer entre deux écueils, le point de vue exprimé par Árpád \citet{sz69,sz00} qui défend l'idée d'une philosophie éléate déterminante pour la constitution des mathématiques grecques, et le point de vue de G. E. L.~\citet{owen57}, qui le réfute.

Plusieurs historiens ont voulu retracer les circonstances de
l'élaboration des \emph{Éléments} d'Euclide. Parmi ceux-ci, Szabó se démarque en considérant la structure axiomatique
des \emph{Éléments} comme une réaction à la philosophie de Parménide.

L'impact de la philosophie éléate sur les mathématiques grecques a été un sujet d'intenses discussions dont le point culminant sont peut-être les débats de la section~II de la conférence de Pise de 1978 \citep{hintikkagruenderagazzi81}, \emph{The early history of the
axiomatic method}
.

Les traces du monde sensible dans le monument de la science grecque que sont les \emph{Éléments} d'Euclide sont circonscrites aux premières définitions: \og un \emph{point} est ce dont il n'y a aucune partie\fg, \og une \emph{ligne} est une longueur sans largeur\fg, \og les \emph{limites d'une ligne} sont des points\fg, \og une \emph{ligne droite} est celle qui est placée de manière égale par rapport aux points qui sont sur elle\fg, etc.\ et à la septième notion commune du premier livre: \og les choses qui s'ajustent les unes sur les autres sont égales entre elles\footnote{Euclide d'Alexandrie, \emph{Les \'El\'ements}, Paris, Presses universitaires de France, 1990-2001, tome~1, pages~150-154 et~178.}\fg\nocite{vitrac9001}.  La perfection de la géométrie d'Euclide est allée de pair avec le développement d'une philosophie qui exclut le réel de la construction des mathématiques.


Les mathématiques grecques rendent compte dans leur forme d'un mode de présentation dont on est en droit de présumer qu'il est une réponse à un problème épistémologique concret. C'est d'autant plus le cas que  
les définitions des objets de la géométrie 
sont davantage des descriptions et se révèlent inertes pour le discours mathématique.
Ils sont \og là\fg\ et on peut même supposer qu'Euclide a fait un effort particulier pour les décrire de manière statique et non construite.

Les nombres, quant à eux, résultent d'une construction consciente: pour paraphraser les deux
premières définitions du septième livre des \emph{Éléments}, il faut \og dire qu'une chose est une\fg\ pour avoir affaire à une unité; il faut \og former une multitude d'unités\footnote{Ibid., tome~2, page~247.}\fg\ pour avoir affaire à un nombre.

Et pourtant, on peut présumer que la question de ce que c'est qu'être un point, une ligne, etc., a agité les mathématiciens grecs. La réponse opératoire qui consiste à dire que les objets sont définis par leur utilisation, par la manière dont on opère sur eux dans les démonstrations, est certainement un renoncement. 

Les postulats des \emph{Éléments} sont des \og demandes\fg\ qui posent le cadre de la validité des raisonnements et ne signifient pas un engagement sur la nature des objets de la géométrie: on peut davantage les lire comme un avertissement sur les limitations de l'application de la géométrie à l'espace sensible. Les raisonnements eux-mêmes réussissent à se placer à un niveau purement discursif et ne font qu'opérer sur les objets sans interroger leur nature ni les rapports qu'ils entretiennent entre eux.

Les témoignages de ce renoncement sont rares et c'est aussi à ce titre que les paradoxes de Zénon sont précieux.

Il apparait ainsi que la science mathématique toute entière est du domaine des opinions et que
le rapport des mathématiques à la vraisemblance est celui d'un devenir
. La forme axiomatique de la géométrie d'Euclide peut être comprise dans ce sens comme une manière d'organiser les opinions sur le monde sensible de manière vraisemblable sans jamais s'exprimer sur la nature du continu, c'est-à-dire sur le rapport entre l'espace et ses parties
.


\bibliography{parmenide-zenon,references-href,Biblio-histoire-des-maths}

\end{document}